\documentclass[11 pt]{amsart}
\usepackage{amssymb}
\usepackage{amsmath}
\usepackage{amsfonts}
\usepackage{graphicx}
\usepackage{amsthm}

\usepackage{enumerate}
\usepackage[mathscr]{eucal}
\usepackage{verbatim}

 \theoremstyle{plain}
\newtheorem{theorem}{Theorem}[section]
\newtheorem{lemma}[theorem]{Lemma}
\newtheorem{proposition}[theorem]{Proposition}

\numberwithin{equation}{section}
\theoremstyle{plain}

\theoremstyle{remark}

\newcommand {\conv}{\mathop{\mathrm{conv}}}

\def\bbR{{\mathbb {R}}}
\def\bbZ{{\mathbb {Z}}}
\def\bbN{{\mathbb {N}}}

\begin{document}


\title
{Dimensions of sums with self-similar sets}

\begin{abstract}
For some self-similar sets $K\subset\bbR^d$ we obtain certain lower bounds for the lower Minkowski dimension of $K+E$ in terms of the lower Minkowski dimension of $E$.
\end{abstract}

\author[]
{Daniel Oberlin and Richard Oberlin}


\address
{Daniel  Oberlin \\
Department of Mathematics \\ Florida State University \\
 Tallahassee, FL 32306}
\email{oberlin@math.fsu.edu}

\address
{Richard Oberlin \\
Department of Mathematics \\ Florida State University \\
 Tallahassee, FL 32306}
\email{roberlin@math.fsu.edu}


\thanks{D.O. was supported in part by NSF Grant DMS-1160680
and R.O. was supported in part by NSF Grant DMS-1068523.}




\maketitle

\section{Introduction}

Suppose $K$ and $E$ are compact subsets of $\bbR^d$ and consider the sum set 
$K+E \dot=\{k+e:k\in K,\, e\in E\}$. We are interested in what can be said about 
$\dim (K+E)$ in terms of $\dim (K)$ and $\dim (E)$.  
Of course there is the obvious lower bound  
\begin{equation}\label{1}
\dim (K+E)\geq\min\{\dim (K),\dim (E)\}.
\end{equation}
If $d\geq 2$, examples involving lower dimensional subspaces show that equality may hold in \eqref{1}, while the less obvious existence of such examples in $\bbR$ is established (for Hausdorff dimension) in \cite{SS}. 

On the other hand, it is clear that 
\begin{equation}\label{.5}
\dim (K+E)\leq \min\{\dim (K)+\dim (E),d\},
\end{equation}
and it is easy to find trivial examples for which equality holds. More interestingly,
if $K,E\subset \bbR$ are classical Cantor sets with ratios of dissection $r_K$ and $r_E$, and if $\log (r_K )/\log (r_E )$ is irrational, then it is shown in \cite{PS} that equality holds in \eqref{.5}. Also, \cite{O3} contains the easy observation that if $K\subset\bbR^d$ is a Salem set, then, for Hausdorff dimension, there is always equality in \eqref{.5}. Here, as in \cite{O3}, we are interested in focusing on particular sets $K\subset\bbR^d$ and finding lower bounds  
\begin{equation*}
\dim (K+E)\geq \Phi \big(K, \dim (E)\big )
\end{equation*}
%
which improve on \eqref{1}. In this note we will be interested in the case when $K$ is self-similar and $\dim =\dim_m$, the lower Minkowski dimension. In particular, we will show that 
for certain classes of self-similar sets
$K\subset\bbR^d$ there exists
$\gamma =\gamma (K)\in (0,1)$ such that 
\begin{equation}\label{conc}
\dim_m (K+E) \ge \gamma d +(1-\gamma )\dim_m (E).
\end{equation}
Such estimates improve on the trivial estimate \eqref{1} when 
\begin{equation*}
\dim_m (E) > \big(\dim_m (K) -\gamma d \big)/(1-\gamma ).
\end{equation*}
We will also obtain more specific results of this type for certain Cantor-like subsets of $\bbR$. Some of the tools we will use are already present in \cite{O1}, \cite{O2}, and \cite{O3}.

\section{Results}

A {\it similarity} on $\bbR^d$ is a map $\phi: x\mapsto rOx+b$ where $r\in (0,1)$, $O$ is an orthogonal transformation of $\bbR^d$, and $b\in\bbR^d$. We will say that a nonempty compact set $K\subset\bbR^d$ is {\it self-similar} if there are $J\in\bbN$ and similarities $\phi_j (x)=r_j O_j x+b_j$, $0\le j\le J$, such that 
\begin{equation}\label{defK}
K=\cup_{j=0}^J \phi_j (K) .
\end{equation}
Note that if the similarities $\phi_j$ map some affine  hyperplane $\mathcal P$ into itself, then $K$ lies in $\mathcal P$. Thus
\begin{equation*}
\dim_m (K+E)>\dim_m (E) 
\end{equation*}
may fail even when $\dim_m (E) $ is close to $d$.
One way to prevent this is to assume that the convex hull $\text{conv} \{b_0 ,\dots b_J \}$ of  $\{b_0 ,\dots b_J \}$ contains an interior point. With this assumption we have the following result. 
\begin{theorem}\label{theorem1}
With $K$ as above, assume that $\phi_j (x)=rx+b_j$ for some $r\in (0,1)$ and for $b_0 ,\dots ,b_J$ such that
$\text{conv} \{b_0 ,\dots b_J \}$ contains an interior point. Suppose $k$ is a positive integer 
such that
\begin{equation*}
k+1\ge \frac{d}{r}.
\end{equation*}
Then
\begin{equation*}
\dim_m (K+E)\geq \frac{d}{k}+\frac{k-1}{k}\dim_m (E) 
\end{equation*}
for, say, compact $E\subset\bbR^d$.
\end{theorem}
\noindent The proof of Theorem \ref{theorem1} is an immediate consequence of the following three results (which will be proved in \S 3). The first of these is implicit in \cite{O3}.

\begin{lemma}\label{lemma1} Let $m_d$ denote Lebesgue measure on $\bbR^d$. Suppose the compact set $K\subset\bbR^d$ satisfies the following condition for some $c_K >0$, some $\gamma\in (0,1)$, and all open $S\subset\bbR^d$: 
\begin{equation}\label{2.02}
m_d (K+S)\geq c_K \, m_d (S)^{1-\gamma}. 
\end{equation}
Then 
\begin{equation}\label{2.2}
\dim_m (K+E)\geq \gamma d +(1-\gamma )\dim_m (E)
\end{equation}
for, say, compact $E\subset\bbR^d$.
\end{lemma}

\noindent Thus our strategy will be to study certain inequalities of the form \eqref{2.02}. (Such estimates were already the subjects of \cite{O1} and \cite{O2}.) One approach to such inequalities is given by the next lemma.
\begin{lemma}\label{lemma2}
(a) Suppose $G$ is an abelian group and $K,S\subset G$. Then, for $k=1,2,\dots$, we have the inequality 
\begin{equation*} 
|\{(x_0 -x_1 ,x_0 -x_2 ,\dots ,x_0 -x_k ):x_i \in K\}|^{1/(k+1)} |S|^{1-k/(k+1)}\le |K+S|
\end{equation*}
where $|\cdot |$ denotes cardinality.

\noindent (b) (Pl\"unnecke-Rusza estimates) With $K$ and $S$ as in (a), for $k=2,3,\dots$,  and with 
\begin{equation*}
K\pm K\pm \cdots \pm K
\end{equation*}
denoting any one of the $2^{k-1}$ possibilities resulting from using $k$ copies of $K$,
there is the inequality
\begin{equation*} 
|K\pm K\pm \cdots \pm K |^{1/k} |S|^{1-1/k}\le |K+S| .
\end{equation*}

\noindent (c) There is a positive constant $C_d$ such that if $K,S\subset \bbR^d$ with $K$ compact and $S$ open then,
for $k=1,2,\dots$, we have the inequality 
\begin{equation*} 
m_{dk}(\{(x_0 -x_1 ,x_0 -x_2 ,\dots ,x_0 -x_k ):x_i \in K\})^{1/(k+1)} m_d (S)^{1-k/(k+1)}\le C_d \,m_d (K+S).
\end{equation*}

\noindent (d) With $K$ and $S$ as in (c) and $k=2,3,\dots$ we have
\begin{equation}\label{3.4} 
m_d (K\pm K\pm \cdots \pm K)^{1/k} m_d (S)^{1-1/k}\le 
C_d \,m_d (K+S).
\end{equation}
\end{lemma}

\noindent Here are two remarks on Lemma \ref{lemma2}. First, we will use only (c) and (d) of 
Lemma \ref{lemma2} but have included (a) and (b) in the statement of the lemma instead of its proof because we wish to draw attention to the possibility that (a) and (b) are instances of a larger family of interesting additive-combinatorial estimates. For example, with $K$ and $S$ as in (a) and (b), we 
{\bf conjecture} the estimate 
\begin{equation}
|\{ (x_1 +x_2 +x_3 ,x_1 +x_4 +x_5 ):x_i \in K\}|^{1/5}|S|^{1-2/5}\le |K+S|.
\end{equation}
Second, an examination of the proof of Proposition 4 in \cite{O1} yields an alternate proof of (c) with $C_d =1$. It seems reasonable that $C_d =1$ should work in (d) as well, but we have not proved this.

\begin{proposition}\label{lemma3}
Suppose $K\subset \bbR^d$ is compact and nonempty and 
satisfies \eqref{defK} where $\phi_j (x) =rx+b_j$ for $r\in (0,1)$ and
$\text{conv} \{b_0 ,\dots b_J \}$ has nonempty interior. Then if 
$$
k+1\ge \frac{d}{r}
$$
it follows that the $k$-fold sum
$$
k\cdot K\doteq K+\cdots +K
$$
has nonempty interior.
\end{proposition}

\noindent We note that Proposition \ref{lemma3} is a higher-dimensional analog of Corollary 2.3 in \cite{CHM}.

Theorem \ref{theorem1} shows that many self-similar sets satisfy estimates of the form \eqref{conc}. But even in the one-dimensional case we do not know an example in which Theorem \ref{theorem1} yields a sharp result. By way of illustration, we consider the Cantor middle-thirds set $C$. The immediate consequence of Theorem \ref{theorem1} for $C$ is that 
\begin{equation}\label{2.105}
\dim_m (C+E)\ge \frac{1}{ 2}+\frac{ 1}{2}\dim_m (E)
\end{equation}
for compact $E\subset \bbR$. But, as was noted in \cite{O3}, the stronger inequality 
\begin{equation}\label{2.106}
\dim_m (C+E)\ge \frac{\log 2}{\log 3}+\big(1-\frac{\log 2}{\log 3}\big)\dim_m (E)
\end{equation}
follows from Lemma \ref{lemma1} and a result in \cite{O2}. 
We have no reason to suspect that \eqref{2.106} is sharp, though it is the best result that our approach (based on Lemma \ref{lemma1}) can give.

One class of sets which generalize $C$ is the collection of homogeneous Cantor sets 
$C_a$, $0<a<1/2$, where 
\begin{equation*}
C_a =\big\{(1-a)\sum_{j=0}^\infty \epsilon_j a^j :\epsilon_j \in \{0,1\}\big\}.
\end{equation*}
The best result of the form \eqref{conc} that we know for the entire class of sets $C_a$ is the following. 
\begin{theorem}\label {prop1}
If $k$ is a positive integer satisfying $k\ge (1-a)/a$ then 
\begin{equation}\label{2.0001}
\dim_m (C_a +E)\geq \frac{1}{k}+\frac{k-1}{k}\dim_m (E) 
\end{equation}
for compact $E\subset \bbR$.
\end{theorem}
\noindent The proof is a consequence of Lemma \ref{lemma1}, (d) of Lemma \ref{lemma2}, 
and the fact that $k\cdot C_a $ contains a nontrivial interval (which follows from either  Corollary 2.3 in \cite{CHM} or (its generalization)
Proposition \ref{lemma3}.) We note that the conclusion \eqref{2.0001} is strongest when $k=2$, yielding then the analog
\begin{equation*}
\dim_m (C_a +E)\ge \frac{1}{ 2}+\frac{ 1}{2}\dim_m (E)
\end{equation*}
of \eqref{2.105}. This follows from Theorem \ref{prop1} only for $1/3\le a<1/2$. However \eqref{2.106} suggests the 
{\bf conjecture} 
\begin{equation}
\dim_m (C_a +E)\ge \frac{\log 2}{\log 1/a}+\big(1-\frac{\log 2}{\log 1/a}\big)\dim_m (E).
\end{equation}
%


The Cantor sets $C_a$ are the self-similar sets corresponding to the two similarities 
$\phi_0 (x) =ax$, $\phi_1 (x) =ax+(1-a)$. The next result concerns the related collections of similarities 
%
\begin{equation}\label{R1}
\phi_j (x)=rx +j ,\ 0\le j\le J.
\end{equation}
\begin{theorem}\label{R2}
Fix a positive integer $J>2$ and suppose that $$\frac{2}{3J}\le r <\frac{1}{J+1}.$$ With the 
$\phi_j$'s as in \eqref{R1}, suppose that the nonempty compact set $K\subset \bbR$ satisfies $K=\cup_{j=0}^J \phi_j (K)$.  Then
\begin{equation*}
\dim_m (K+E)\ge \frac{2}{3}+\frac{1}{3} \dim_m (E)
\end{equation*}
for compact $E\subset \bbR$.
\end{theorem}
\noindent Theorem \ref{R2} follows from Lemma \ref{lemma1}, (c) of Lemma \ref{lemma2}, and the following result:

\begin{proposition}\label{R3} With $K$ as in Theorem \ref{R2}, the subset 
$$
\{(x_0 -x_1 ,x_0 -x_2 ): x_i \in K\}
$$
of $\bbR^2$ has nonempty interior.
\end{proposition}

Another class of Cantor-like sets is obtained as follows: fix a positive integer $n\ge 3$ and a subset $A$ of $\{0,1,\dots ,n-1\}$. Define 
\begin{equation*}
C_{n,A}=\{\sum_{j=1}^\infty \frac{a_j}{n^j}:a_j \in A\}.
\end{equation*}
The following generalization of \eqref{2.106} was proved in \cite{O3}: if $0\in A$ and $|A|=n-1$ then 
\begin{equation*}
\dim_m (C_{n,A} +E)\ge \frac{\log (n-1)}{\log n}+\big(1-\frac{\log (n-1)}{\log n}\big)\dim_m (E).
\end{equation*}
Here is another result for the sets $C_{n,A}$ (without the restriction $|A|=n-1$).
\begin{theorem}\label{theorem2}
Let $G_n$ stand for the group of integers modulo n. 
Suppose that $A\subset G(n)$.
For fixed $k=1,2,\dots$, suppose that 
\begin{equation}\label{4.2}
\{(a_0 -a_1 , a_0 -a_2 ,\dots ,a_0 -a_k ): a_j \in A\}=(G_n )^k .
\end{equation}
Then 
\begin{equation*}
\dim_m (C_{n,A} +E)\geq  \frac{k}{k+1} +\frac{\dim_m (E)}{k+1}.
\end{equation*}
\end{theorem}
\noindent The proof is a direct consequence of (c) of Lemma \ref{lemma2}, Lemma \ref{lemma1}, and the following result, to be proved in \S 3.
\begin{proposition}\label{proposition2}
If \eqref{4.2} holds then 
\begin{equation}\label{2.222}
m_k \big(\{(x_0 -x_1 ,x_0 -x_2 ,\dots ,x_0 -x_k ):x_j \in C_{n,A} \}\big)>0 .
\end{equation}
\end{proposition}
There is an alternative approach to Theorem \ref{theorem2} based on (a) of Lemma \ref{lemma2} and Theorem 2 of \cite{O2}. We choose to prove Theorem \ref{theorem2} based on Proposition \ref{proposition2} in order to establish (alone with Proposition \ref{R3}) some motivation for our {\bf conjecture} that, for given $k$, we have
\begin{equation}\label{2.22}
m_k \big(\{(x_0 -x_1 ,x_0 -x_2 ,\dots ,x_0 -x_k ):x_j \in C_a \}\big)>0 
\end{equation}
 so long as the parameter $a$ is close enough to $1/2$. If true, this conjecture might 
 provide, via (c) of Lemma \ref{lemma2} and Lemma \ref{lemma1},  
 an improvement on Theorem \ref{prop1}.

\section{Proofs}

{\it Proof of Lemma \ref{lemma1}:}
For $E\subset\bbR^d$, the condition $\dim_m (E)\geq \beta$ is equivalent to the estimate 
\begin{equation}\label{2.3}
m_d \big(E+B(0,\delta ) \big)\gtrsim \delta^{d-\beta+\epsilon}
\end{equation}
where the implied constant depends on $\epsilon >0$. 
Inequalities \eqref{2.3} and \eqref{2.02} together imply
\begin{equation*}
m_d \big(K+E+B(0,\delta ) \big)\gtrsim \delta^{(d-\beta+\epsilon )(1-\gamma)} =
\delta^{d-\big(\beta +\gamma (d-\beta)\big)+\epsilon (1-\gamma)}
\end{equation*}
which, upon replacing $E$ by $K+E$ in \eqref{2.3}, yields \eqref{2.2}.

$\ \ \ \ \ \ \ \ \ \ \ \ \ \ \ \ \ \ \ \ \ \ \ \ \ \ \ \ \ \ \ \ \ \ \ \ \ \ \ \ \ \ \ \ \ \ \ \ \ \ \ 
\ \ \ \ \ \ \ \ \ \ \ \ \ \ \ \ \ \ \ \ \ \ \ \ \ \ \ \ \ \ \ \ \ \ \ \ \ \ \ \ \ \ \ \ \ \ \ \ \ \ \
\square$

{\it Proof of Lemma \ref{lemma2}:}
To see (a) we assume that $K$ is finite and let
\begin{equation*}
\{(x_0^n, x_1^n ,\dots ,x_k^n ):1\le n\le N,\, x_i^n \in K\}
\end{equation*}
be such that
\begin{equation*} 
\{(x_0^n -x_1^n ,\dots ,x_0^n -x_k^n ):1\le n\le N\}
\end{equation*}
is a one-to-one enumeration of 
\begin{equation*}
\{(x_0 -x_1 ,x_0 -x_2 ,\dots ,x_0 -x_k ):x_i \in K\}. 
\end{equation*}
Then one easily checks that the map
\begin{equation*}
(n,s)\mapsto (x_0^n +s ,\dots ,x_k^n +s) 
\end{equation*}
is a one-to-one mapping of $\{1,2,\dots ,N\}\times S$ into $(K+S)^{k+1}$.

Part (b) is just a restatement of the Pl\"unnecke-Rusza estimates (Corollary 6.29 in \cite{TV}) which say that if $C$ is any positive constant satisfying $|K+S|\le C\,|S|$ then 
\begin{equation*}
|K\pm K\pm \cdots \pm K |\le C^k |S|.
\end{equation*}

Next we will give the proof for (d) - part (c) can be proved similarly (but see also the remarks immediately following this proof). The proof is just an approximation argument based on (b), but we include it because it is not completely straightforward. Let $\mathcal L _n$ be the additive group in $\bbR^d$ generated by the scaled unit vectors $(1/n) \,u_j , 1\le j\le d$. If $E\subset\bbR^d$ is, for example,  a finite union of rectangles 
$
\prod [a_j ,b_j ]
$,
 then
\begin{equation}\label{3.45}
m_d (E)=\lim_{n\rightarrow \infty}\frac{1}{n^d} |\mathcal L _n \cap E|.
\end{equation}
Suppose to begin that $K$ and $S$ are finite unions of such closed 
and nondegenerate
rectangles. If $x\in K\pm K\pm \cdots \pm K$ then, for large $n$,  there are 
$\ell_1 ,\dots ,\ell_k \in \mathcal L _n \cap K$ such that $|x-(\ell_1 \pm \ell_2 \pm \cdots \pm \ell_k )|\le c_d /n$. Thus, for some $C_d \ge 1$, 
\begin{equation*}
m_d (K\pm K\pm \cdots \pm K)\le \frac{C_d}{n^d}\, |(\mathcal L _n \cap K )\pm
(\mathcal  L _n \cap K )\pm\cdots \pm (\mathcal  L _n \cap K )|.
\end{equation*}
Therefore
\begin{multline*} 
m_d (K\pm K\pm \cdots \pm K)^{1/k}\Big(\frac{|\mathcal L_n \cap S |}{n^d}\big)^{1-1/k}\le \\
\Big( \frac{C_d}{n^d}\, |(\mathcal L _n \cap K )\pm
(\mathcal  L _n \cap K )\pm\cdots \pm (\mathcal  L _n \cap K)|\Big)^{1/k}\Big(\frac{|\mathcal L_n \cap S |}{n^d}\big)^{1-1/k} \le \\
\frac{C_d}{n^d} |( \mathcal L _n \cap K) +(\mathcal L_n \cap S)|\le \frac{C_d}{n^d} | 
\mathcal L _n \cap (K+S)|, 
\end{multline*}
where (b) was used to obtain the next-to-last inequality. Letting $n\rightarrow\infty$ and using  \eqref{3.45} gives \eqref{3.4} when $K$ and $S$ are finite unions of closed  rectangles. If $K\subset \bbR^d$ is compact, then $K=\cap_j K_j$ where 
$K_1 \supset  K_2 \supset \cdots$ and the $K_j$'s are finite unions of closed rectangles. 
With $S$ still a finite union of closed rectangles, suppose that $G$ is open and $K+S\subset G$. Then $K_j +S\subset G$ for large $j$ and so
\begin{multline*} 
m_d (K\pm K\pm \cdots \pm K)^{1/k}m_d (S)^{1-1/k}\le \\
m_d (K_n \pm K_n \pm \cdots \pm K_n )^{1/k}m_d (S)^{1-1/k}\le \\
C_d \,m_d (K_j +S) \leq C_d \,m_d (G)
\end{multline*}
for large $j$. Taking an infimum over open $G$ with $K+S\subset G$ shows that
\eqref{3.4} holds whenever $K$ is compact and $S$ is a finite union of closed rectangles. Approximating open rectangles from inside by closed rectangles gives the result when $S$ is a finite union of open rectangles and then \eqref{3.4} follows whenever $K$ is compact and $S$ is open.

$\ \ \ \ \ \ \ \ \ \ \ \ \ \ \ \ \ \ \ \ \ \ \ \ \ \ \ \ \ \ \ \ \ \ \ \ \ \ \ \ \ \ \ \ \ \ \ \ \ \ \ 
\ \ \ \ \ \ \ \ \ \ \ \ \ \ \ \ \ \ \ \ \ \ \ \ \ \ \ \ \ \ \ \ \ \ \ \ \ \ \ \ \ \ \ \ \ \ \ \ \ \ \
\square$

{\it Proof of Proposition \ref{lemma3}:}
A translation argument based on the observation 
\begin{equation*}
K+b =\cup_{j=0}^J \big(r(K+b) +b_j +(1-r )b\big) 
\end{equation*}
shows that we can assume that $b_0 =0$. By our assumption that
$\text{conv} \{b_0 ,\dots b_J \}$ has nonempty interior
we can relabel to ensure that 
$\{b_1 ,\dots ,b_d \}$ is a maximal linearly independent subset of $\{b_1 ,\dots ,b_J \}$.
For any $F\subset \bbR^d$ we will write 
\begin{equation*}
T(F)=\cup_{j=0}^J (rF+b_j ).
\end{equation*}
It is well known that if $K_0$ is any compact set satisfying 
$$
rK_0 +b_j
\subset K_0
$$
for $j=0,\dots ,J$  and if 
\begin{equation*}
K_{n+1}=T(K_n )
\end{equation*}
then $\cap_n K_n =K$. Fix such a $K_0$ with 
$$
\text{conv} \{(k+1)b_0 ,\dots (k+1)b_d \}
\subset k\cdot K_0 
$$
(a large enough closed ball with center at the origin will suffice) .
Now assume that $n\ge 0$ and
\begin{equation}\label{N+1}
\text{conv} \{(k+1)b_0 ,\dots (k+1)b_d \}
\subset k\cdot K_n . 
\end{equation}
We will show that then 
\begin{equation}\label{N+2}
\text{conv} \{(k+1)b_0 ,\dots (k+1)b_d \}
\subset k\cdot T(K_n )=k\cdot K_{n+1} .
\end{equation}
This implies, via an easy compactness argument, that  
\begin{equation*}
\text{conv} \{(k+1)b_0 ,\dots (k+1)b_d \}
\subset \cap_{n=0}^\infty (k\cdot K_n )=k\cdot K ,
\end{equation*}
yielding the conclusion of Proposition \ref{lemma3}. 

It remains to show that \eqref{N+1} implies \eqref{N+2}.
For $\theta >0$ define 
$$
\text{sfloor}(\theta)=\max\{p\in\bbZ :p<\theta \}
$$
and $\text{sfloor}(0)=0$. If $x\in \text{conv} \{(k+1)b_0 ,\dots (k+1)b_d \}$ then 
we can write, for $\theta_j \ge 0$, $\sum_1^d \theta_j \leq k+1$,  
\begin{equation*}
x=\sum_{j=1}^d \theta_j b_j =
\sum_{j=1}^d \big(\theta_j -\text{sfloor}(\theta_j )\big)b_j +
\sum_{j=1}^d \text{sfloor}(\theta_j )b_j .
\end{equation*}
We note that 
\begin{equation*}
\sum_{j=1}^d \big(\theta_j -\text{sfloor}(\theta_j )\big)\le d,\ \ 
\sum_{j=1}^d \text{sfloor}(\theta_j )\le k .
\end{equation*}
Thus, recalling that $k\cdot \{b_0 ,\dots ,b_J\}$ indicates the $k$-fold sum
of $ \{b_0 ,\dots ,b_J\}$, that $d\le (k+1)r$, and that $b_0 =0$, we see that 
\begin{multline*}
x\in d\ \text{conv} \{ b_0 ,\dots  b_d \}+k\cdot \{b_0 ,\dots ,b_d\}\subset \\
r\ \text{conv} \{(k+1)b_0 ,\dots (k+1)b_d \}+k\cdot \{b_0 ,\dots ,b_J\} \subset \\
r(k\cdot K_n )+k\cdot \{b_0 ,\dots ,b_J\}=k\cdot T(K_n ),
\end{multline*}
where we have used \eqref{N+1} for the last inclusion. This yields \eqref{N+2} and therefore concludes the proof of the Proposition \ref{lemma3}.

$\ \ \ \ \ \ \ \ \ \ \ \ \ \ \ \ \ \ \ \ \ \ \ \ \ \ \ \ \ \ \ \ \ \ \ \ \ \ \ \ \ \ \ \ \ \ \ \ \ \ \ 
\ \ \ \ \ \ \ \ \ \ \ \ \ \ \ \ \ \ \ \ \ \ \ \ \ \ \ \ \ \ \ \ \ \ \ \ \ \ \ \ \ \ \ \ \ \ \ \ \ \ \
\square$

{\it Proof of Proposition \ref{R3}:}
Beginning with some notation, for $F\subset\bbR$ we define
\begin{multline*}
T(F)\dot=\bigcup_{j=0}^J \phi_j (F)=\bigcup_{j=0}^J \big(rF+j),\\ 
D(F)\dot=\{(x_0 -x_1 ,x_0 -x_2 ): x_i \in F \}.
\end{multline*}
For the moment we believe the following result:
\begin{lemma}\label{R4}
Suppose $r \geq \frac{2}{3J}$, and write 
\[
V = \{(0,0), (0,J), (J,0), (J,2J), (2J,J), (2J, 2J)\}.
\]
Suppose for some $w \in \bbR^2$ for some $F\subset\bbR$ that $D(F)$ contains $w+\conv(V).$ Then $D(T(F))$ contains $rw + (1-J,1-J) + \conv(V).$
\end{lemma}

\noindent Let $K_0$ be a closed interval containing $0$ and large enough that 
\[
(\frac{1-J}{1-r},\frac{1-J}{1-r}) + \conv(V) \subset D(K_0).
\]
With $K_{n+1} = T(K_n)$, so that $K = \bigcap_{n}K_n$,
the lemma yields
\[
(\frac{1-J}{1-r},\frac{1-J}{1-r}) + \conv(V) \subset D(K_n)
\]
for each $n$. Thus, by a compactness argument, 
\[
(\frac{1-J}{1-r},\frac{1-J}{1-r}) + \conv(V) \subset D(K).
\]
Since $\conv(V)$ has nonempty interior, the conclusion of Proposition \ref{R3} follows.

$\ \ \ \ \ \ \ \ \ \ \ \ \ \ \ \ \ \ \ \ \ \ \ \ \ \ \ \ \ \ \ \ \ \ \ \ \ \ \ \ \ \ \ \ \ \ \ \ \ \ \ 
\ \ \ \ \ \ \ \ \ \ \ \ \ \ \ \ \ \ \ \ \ \ \ \ \ \ \ \ \ \ \ \ \ \ \ \ \ \ \ \ \ \ \ \ \ \ \ \ \ \ \
\square$

{\it Proof of Lemma \ref{R4}:}

Write 
\[
P(\{0, \ldots, J\}) = \{(x_0 + x_1, x_0 + x_2) : x_i \in \{0, \ldots, J\}, i = 0 ,1, 2\}.
\]
Observe that
\begin{align*}
D(T(F)) &= rD(F) + D(\{0, \ldots, J\}) \\ 
&=rD(F) + P(\{0, \ldots, J\}) + (-J,-J) \\
&\supset rw + r\conv(V) + P(\{0, \ldots, J\}) + (-J,-J).
\end{align*}
Then Lemma \ref{R4} will follow when we show that
\begin{equation} \label{wcontains}
W\dot= \,r\conv(V) + P(\{0, \ldots, J\}) \supset (1,1) + \conv(V).
\end{equation} 

To prove \eqref{wcontains} it is sufficient (and necessary) to establish the inclusions
\begin{align} \label{firstjtileunion}
& (1,1) + \conv(\{(0,0), (J,0),(0,J),(J,2J)\})\\
\nonumber &=  \conv(\{(1,1),(J+1,1),(1,J+1),(J+1, 2J+1)\}) \subset W
\end{align}
and 
\begin{align} \label{secondjtileunion}
&(1,1) + \conv(\{(J,0), (J,2J),(2J,J),(2J,2J)\}) \\
\nonumber &=  \conv(\{(J+1, 1),(J+1,2J+1),(2J+1,J+1),(2J+1, 2J+1)\}) \\
\nonumber&\subset W .
\end{align}
To this end we define the triangles 
\begin{equation*}
A_u \dot= \conv\big(\{(0,0),(1,1),(0,1)\}\big),\  A_l \dot= \conv\big(\{(0,0),(1,1),(1,0)\}\big).
\end{equation*}
To see \eqref{firstjtileunion} and \eqref{secondjtileunion} we will need the inclusions 
\begin{equation} \label{utcont}
A_u \subset r\conv(V) + \{(0,0), (-1,0)\}
\end{equation}
and 
\begin{equation} \label{ltcont}
A_l \subset r\conv(V) + \{(0,0), (0,-1)\}.
\end{equation} 
Here is the proof of \eqref{utcont} (the proof of \eqref{ltcont} follows from a symmetric argument). 

Our assumption $r\ge \frac{1}{2J}$ implies that 
\begin{equation*}
\conv\big(\{(0,0), (0,rJ), (1-rJ,1), (1,1)\}\big) \subset r\conv (V).
\end{equation*}
Thus it suffices to show that 
\begin{equation*} \label{topgap}
\conv(\{(0,1),(0,rJ),(1-rJ,1)\} \subset r\conv(V) + (-1,0).
\end{equation*}
By convexity this will follow from
\[
\big\{\big(\frac{1}{r},\frac{1}{r}\big),\big(\frac{1}{r},J\big), \big(\frac{2}{r} - J,\frac{1}{r}\big)
\big\} \subset \conv(V).
\]
The points $(\frac{1}{r},\frac{1}{r}),(\frac{1}{r},J)$ can be checked by using $r \geq \frac{1}{2J}$. For the remaining point, we use  $\frac{2}{3J}\leq r \leq \frac{1}{J}$ to see that $J \leq \frac{2}{r} - J \leq 2J$. Then the condition $r \geq \frac{1}{2J}$ ensures that $(\frac{2}{r} - J,\frac{1}{r})$ lies above the line of slope 1 passing through $(J,0)$. 
This establishes \eqref{utcont}.
(Here, and in the remainder of this proof, a picture may be helpful.)

To apply \eqref{utcont} and \eqref{ltcont} we first observe that
for each integer $j$ with $1 \leq j \leq J$, for each $k$ such that $1 \leq k \leq J-1 + j$, and with $l = \min(j,k)$ we have
\begin{align*}
&\{(j,k), (j-1,k), (j,k-1)\} \\
&=\{(l + (j-l), l + (k - l)), ((l-1) + (j-l),(l-1) + (k+1-l)), \\
&((l-1) + (j+1-l),(l-1) + (k-l))\}\\
&\subset P(\{0, \ldots, J\})
\end{align*}
and
\[
\{(j,J+j),(j,J-1+j)\} \subset P(\{0, \ldots, J\}).
\]
It then follows from \eqref{utcont} and \eqref{ltcont} that 
$$
(A_u \cup A_l )+(j,k)\subset r\conv (V) +P(\{0, \ldots, J\})=W
$$
and similarly that $A_l +(j,J-j)\subset W$. 
Unioning over $j$, $1\le j\le J$, gives
\begin{equation*}
 \conv(\{(1,1),(J+1,1),(1,J+1),(J+1, 2J+1)\}) \subset W
\end{equation*}
which is \eqref{firstjtileunion}. 

The proof of  \eqref{secondjtileunion} is similar, starting with the observation that
for $J+1 \leq j \leq 2J$, $j - J + 1 \leq k \leq 2J$, and $l = \max(j-J,k-J)$, we have
\begin{align*}
&\{(j,k), (j-1,k), (j,k-1)\} \\
&=\{(l + (j-l),l+(k-l)), (l + (j-1-l), l + (k-l)),\\
&(l + (j-l),l + (k-1-l))\}\\
&\subset P(\{0, \ldots, J\})
\end{align*}
and
\[
\{(j,j-J),(j-1,j-J)\} \subset P(\{0, \ldots, J\}).
\]
This completes the proof of Lemma \ref{R4}.

$\ \ \ \ \ \ \ \ \ \ \ \ \ \ \ \ \ \ \ \ \ \ \ \ \ \ \ \ \ \ \ \ \ \ \ \ \ \ \ \ \ \ \ \ \ \ \ \ \ \ \ 
\ \ \ \ \ \ \ \ \ \ \ \ \ \ \ \ \ \ \ \ \ \ \ \ \ \ \ \ \ \ \ \ \ \ \ \ \ \ \ \ \ \ \ \ \ \ \ \ \ \ \
\square$

{\it Proof of Proposition \ref{proposition2}:}

For $y\in [0,1]$ we will write $y=.y_1 y_2 y_3 \dots$ if 
$$
y=\sum_{j=1}^\infty y_j n^{-j},\, y_j \in \{0,1,\dots ,n-1\}.
$$
Fix $y_l =.y_{l1}y_{l2}y_{l3}\dots$ for $l=1,2,\dots,k$. Suppose that for some $m\ge 2$ there are $x_0 ,x_1 ,\dots ,x_k \in C_{n,A}$ and $\delta_1 ,\dots , \delta_k \in\{-1 ,0\}$ 
(depending on $m$) with 
\begin{multline*}
(x_0 -x_1 ,x_0 -x_2 ,\dots ,x_0 -x_k ) =\\
(\delta_1 ,\delta _2 ,\dots ,\delta _k) + (.y_{1m}y_{1(m+1)}\dots ,\
.y_{2m}y_{2(m+1)} \dots , \
\dots ,\
.y_{km}y_{k(m+1)}\dots ).
\end{multline*}
Our immediate goal is to show that the same thing is true if $m$ is replaced by $m-1$.
That is, we want to show that 
there are $x'_0 ,x'_1 ,\dots ,x'_k \in C_{n,A}$ and $\delta'_1 ,\dots , \delta'_k \in\{-1 ,0\}$
such that
\begin{multline*} 
(x'_0 -x'_1 ,x'_0 -x'_2 ,\dots ,x'_0 -x'_k ) =\\
(\delta'_1 ,\delta' _2 ,\dots ,\delta' _k) + (.y_{1(m-1)}y_{1m}\dots ,\
.y_{2(m-1)}y_{2m} \dots , \
\dots ,\
.y_{k(m-1)}y_{km}\dots ).
\end{multline*}
By our hypothesis \eqref{4.2} we can choose 
$$
a_0 ,a_1 ,\dots ,a_k \in A
$$ 
such that $a_0 -a_l \equiv y_{l(m-1)} -\delta_l \,\mod n$ for $1\le l \le k$. Then
if 
$$
x' _{p} =(a_p +x_{p} )/n ,\, p=0,1, ,\dots ,k,
$$
we have $x'_p \in  C_{n,A}$
and, for $1\le l\le k$,
\begin{equation*}
x'_0 -x'_l =(a_0 -a_l )/n +\delta_l /n +.0y_{lm}y_{l(m+1)}\dots .
\end{equation*}
Now $a_0 -a_l =y_{l(m-1)} -\delta_l +n\delta'_l$ for some integer $\delta'_l$ and so 
\begin{multline*} 
(x'_0 -x'_1 ,x'_0 -x'_2 ,\dots ,x'_0 -x'_k ) =\\
(\delta'_1 ,\delta' _2 ,\dots ,\delta' _k) + (.y_{1(m-1)}y_{1m}\dots ,\
.y_{2(m-1)}y_{2m} \dots , \
\dots ,\
.y_{k(m-1)}y_{km}\dots ).
\end{multline*}
Because $x'_p \in [0,1]$ for $0\le p\le k$ it follows that $\delta'_l \in\{-1,0\}$ for $1\le l\le k$.

It now follows by induction that if, for each $1\le l\le d$, 
$y_l =.y_{l1}y_{l2}y_{l3}\dots$ where only finitely many of the $y_{lj}$ are nonzero,
then there are $x_0 ,x_1 ,\dots ,x_k \in C_{n,A}$ and $\delta_1 ,\dots , \delta_k \in\{-1 ,0\}$ 
with 
\begin{multline*}
(x_0 -x_1 ,x_0 -x_2 ,\dots ,x_0 -x_k ) =\\
(\delta_1 ,\delta _2 ,\dots ,\delta _k) + (.y_{11}y_{12}\dots ,\
.y_{21}y_{22} \dots , \
\dots ,\
.y_{k1}y_{k2}\dots ).
\end{multline*}
Then a compactness argument shows that the \lq\lq finitely many nonzero" restriction may be removed. Thus $[0,1]^k$ may be covered by $2^k$ translates of 
\begin{equation*}
\{(x_0 -x_1 ,x_0 -x_2 ,\dots ,x_0 -x_k ):x_j \in C_{n,A}\}
\end{equation*}
and \eqref{2.222} follows.

$\ \ \ \ \ \ \ \ \ \ \ \ \ \ \ \ \ \ \ \ \ \ \ \ \ \ \ \ \ \ \ \ \ \ \ \ \ \ \ \ \ \ \ \ \ \ \ \ \ \ \ 
\ \ \ \ \ \ \ \ \ \ \ \ \ \ \ \ \ \ \ \ \ \ \ \ \ \ \ \ \ \ \ \ \ \ \ \ \ \ \ \ \ \ \ \ \ \ \ \ \ \ \
\square$

\end{document}